\documentclass{article}
\usepackage{graphicx} 
\usepackage{xcolor}
\usepackage{amsmath,amssymb,amsthm}
\usepackage{float}
\usepackage{array}
\usepackage{longtable}
\usepackage{booktabs}
\usepackage{hyperref}
\usepackage{graphicx}
\usepackage{eucal}
\usepackage{tikz}
\usepackage{mathtools}
\usepackage{setspace}
\graphicspath{ {images/} }

\newtheorem{theorem}{Theorem}[section]
\newtheorem{lemma}[theorem]{Lemma}
\newtheorem{corollary}[theorem]{Corollary}

\newtheorem{problem}{Problem}

\theoremstyle{definition}
\newtheorem{definition}[theorem]{Definition}

\title{Constructive and computational methods for extremal and near-extremal box colorings}
\author{Jeremy F.~Alm\footnote{\texttt{jeremy.alm@siu.edu}}  ~and Alex Miller\\ School of Mathematical and Statistical Sciences\\
Southern Illinois University\\Carbondale, Illinois, USA}
\date{\today}

\begin{document}

\maketitle

\begin{abstract}
    We investigate extremal and near-extremal 2-colorings of three-dimensional grids without monochromatic axis-parallel boxes. Our approach focuses on the parameter $\rho(a,b)$, defined to be the minimum number of monochromatic rectangles in a 2-coloring of an $a\times b$ grid. We determine $\rho(a,b)$ for several small grids and use these values to obtain upper bounds on the number of layers that can be stacked to form an $a\times b\times c$ grid without creating a monochromatic box. Explicit constructions and computational certificates show that these bounds are attained for several cases, including grids of dimensions $5\times5\times100$, $5\times6\times75$, and $6\times6\times75$. For the grid of dimension  $5\times7\times c$, we construct a box-free coloring with $c=50$, while our counting bound gives $c\leq 52$. This gap suggests that the counting bound may not always be sharp.

We also investigate the largest $n\times n\times n$ cube admitting a box-free 2-coloring. We found box-free colorings of the $n$-cube for $n\leq 10$, and colorings of the $11\times11\times11$ cube with only eight monochromatic boxes. On the other hand, we show that every 2-coloring of the $21\times21\times21$ cube contains at least 21 monochromatic boxes. Thus the largest box-free cube has side length between 10 and 20, inclusive.  We conclude with conjectures and open problems.
\end{abstract}

\section{Introduction}

Let $[n]$ denote $\{1,2,\ldots,n\}$. A three-dimensional grid is the product $[a]\times[b]\times[c]$. A box is a set

$$
\{a_1,a_2\}\times\{b_1,b_2\}\times\{c_1,c_2\}\subseteq [a]\times[b]\times[c].
$$

A $2$-coloring of the grid is a function

$$
f:[a]\times[b]\times[c]\to\{\text{red},\text{blue}\}.
$$

We study extremal and near-extremal $2$-colorings of three-dimensional grids with respect to monochromatic boxes.

The problem has its roots in the classical Zarankiewicz problem \cite{zarankiewicz}. Zarankiewicz asked how many points may be selected from a grid without containing the four corners of an axis-parallel rectangle. Kövári, Sós, and Turán \cite{Kovari} subsequently obtained bounds for this problem, including the asymptotic order of the maximum size of a rectangle-free subset of an $n\times n$ grid, which is $n^{3/2}$. Alm and Manske \cite{Alm} later revisited these results, giving alternative proofs based in part on finite projective planes and recovering the asymptotic result. 

A related line of work studies monochromatic \emph{squares} in colored grids. Axenovich and Manske \cite{AxMan} considered monochromatic squares in $2$-colored  square grids and draw a connection to Van der Waerden numbers. Bacher and Eliahou \cite{Bacher} used computational methods to prove that every $2$-coloring of $[15]\times[15]$ contains such a square, while $[14]\times[14]$ admits a coloring avoiding one. More recently, Apon, Gasarch, and Lawler \cite{complexity} investigated the computational complexity of rectangle-free grid coloring, showing, among other results, that extending a partial coloring to a rectangle-free coloring is NP-complete.

Most directly relevant to the present work are the papers of Cooper, Fenner, and Purewal \cite{Fenner} and of Fenner, Gasarch, Glover, and Purewal \cite{fenner2012rectanglefreecoloringgrids}. Fenner et al. studied the two-dimensional problem of determining exactly which grids $[n]\times[m]$ admit rectangle-free $c$-colorings, obtaining complete classifications for $c=2,3,$ and $4$. Cooper, Fenner, and Purewal considered the higher-dimensional analogue: a grid

$$
[a_1]\times\cdots\times[a_d]
$$

is $c$-guaranteed if every $c$-coloring contains a monochromatic box. They obtained asymptotically tight conditions for a grid to be $c$-guaranteed and studied minimally $c$-guaranteed grids.

Our point of view is complementary. Rather than asking only when a monochromatic box is forced, we ask how close one can come to avoiding boxes and what extremal and near-extremal colorings look like. For several grids we give explicit constructions that give lower bounds that match known upper bounds from \cite{Fenner}. We then supplement these constructions with computational results, which yield further exact results and allow us to explore the minimum number of monochromatic boxes in cubes. In particular, for $n\times n\times n$ cubes we study the minimum number $m(n)$ of unavoidable monochromatic boxes and use min-conflicts search to obtain new bounds.  We ask what is the largest $n$ such that the $n\times n\times n$ cube can be colored without monochromatic boxes and show that $10\leq n < 21$. 

Throughout, we use the minimum number of monochromatic rectangles in two dimensions to investigate the existence and structure of box-free colorings in three dimensions.

\section{Upper Bounds from Cooper et al.}

\begin{definition}
    Let $N(a,b) = \max(\{c : [a]\times [b]\times [c] \text{ admits a box-free $2$-coloring} \})$.
\end{definition}

Then we have the following upper bounds on $c$ for given $a,b$ pairs from \cite{Fenner}. As we will see in Section \ref{sec:layer}, Table \ref{tab:UB}, some of these entries can be improved. 

\begin{table}[H]
\centering
\small
\setlength{\tabcolsep}{4.5pt}
\renewcommand{\arraystretch}{1.08}
\begin{tabular}{c|*{10}{c}}
  & 3 & 4 & 5 & 6 & 7 & 8 & 9 & 10 & 11 & 12 \\
\hline
3  &   &   &     &     & 126 & 84 & 72 & 67 & 66 & 66 \\
\hline
4  &   &   &     &     & 126 & 84 & 72 & 67 & 66 & 66 \\
\hline
5  &   &   & 100 & 75 & 52 & 46 & 45 & 45 & 39 & 36 \\
\hline
6  &   &   & 75  & 75 & 52 & 46 & 45 & 45 & 39 & 36 \\
\hline
7  & 126 & 126 & 52 & 52 & 52 & 45 & 39 & 36 & 33 & 32 \\
\hline
8  & 84 & 84 & 46 & 46 & 45 & 44 & 39 & 36 & 33 & 32 \\
\hline
9  & 72 & 72 & 45 & 45 & 39 & 39 & 36 & 33 & 30 & 29 \\
\hline
10 & 67 & 67 & 45 & 45 & 36 & 36 & 33 & 32 & 30 & 29 \\
\hline
11 & 66 & 66 & 39 & 39 & 33 & 33 & 30 & 30 & 29 & 27 \\
\hline
12 & 66 & 66 & 36 & 36 & 32 & 32 & 29 & 29 & 27 & 27
\end{tabular}
\caption{Upper bounds on $N(a,b)$, the largest $c$
 for which
$[a]\times[b]\times[c]$ admits a box-free $2$-coloring.  These follow
from the $2$-guaranteed bounds of \cite{Fenner} et al.\ by subtracting one.}
\label{tab:fenner-upper-bounds}
\end{table}

\section{The Rectangle Parameter $\rho(a,b)$}

\begin{definition}

    Let $\rho(a,b)$ be the minimum possible number of monochromatic rectangles (counting both colors) in any $2$-coloring of $[a]\times[b]$. 
\end{definition}

The lower bounds in Table 2 are obtained using convexity arguments, sometimes supplemented by an examination of the equality conditions. We illustrate the method below for \(\rho(7,7)\). The matching upper bounds are established by explicit colorings given in Appendix A.  See the data in Table \ref{tab:rho}, which appear in the On-Line Encyclopedia of Integer Sequences \cite{oeis} as A399621. The diagonal sequence $\rho(n,n)$, beginning $0,0,0,0,2,6,21,40$, is recorded separately as A399620.

\begin{table}[hbt]
    \centering
    \begin{tabular}{c|c|c|c|c|c|c}
        $\rho(a,b)$ & 3 & 4 & 5 & 6 & 7 & 8 \\
        \hline 
       3  & 0 & 0 & 0 & 0 & 1 & 2 \\ 
        \hline 
       4  &  & 0 & 0 & 0 & 2 & 4 \\ 
        \hline 
       5  &  &  & 2 & 4 & 8 & 12 \\ 
        \hline 
       6  &  &  &  & 6 & 12 & 18 \\ 
        \hline 
       7  &  &  &  &  & 21 & 30 \\ 
        \hline 
       8  &  &  &  &  &  &  40 \\ 
    \end{tabular}
    \caption{Exact values for $\rho(a,b)$}
    \label{tab:rho}
\end{table}

To illustrate the method, we show that $\rho(7,7)\ge 21$. For each color $c$ and pair of columns $i<j$, let $k_{c,i,j}$ be the number of rows in which both entries have color $c$. The total number of monochromatic rectangles is then
\[
R=\sum_{c,i<j}\binom{k_{c,i,j}}2.
\]

A row containing $r$ red entries and $7-r$ blue entries contributes
\[
\binom r2+\binom{7-r}2\ge 9
\]

to the sum of the $k_{c,i,j}$. Since there are seven rows,
\[
\sum_{c,i<j}k_{c,i,j}\ge 63.
\]

For every nonnegative integer $k$,
\[
\binom{k}{2}
=k-1+\frac12(k-1)(k-2)
\ge k-1,
\]

since $(k-1)(k-2)\ge0$ for every nonnegative integer $k$.

Consequently,

\begin{align*}
R
&=\sum_{c,i<j}\binom{k_{c,i,j}}2\\
&=\sum_{c,i<j}(k_{c,i,j}-1)
+\frac12\sum_{c,i<j}(k_{c,i,j}-1)(k_{c,i,j}-2)\\
&\ge \sum_{c,i<j}k_{c,i,j}-42\\
&\ge63-42=21.
\end{align*}

Thus $\rho(7,7)\ge21$, as claimed.

\section{A Layer Bound for $N(a,b)$}\label{sec:layer}

We get the following upper bound on the number of layers of an $a\times b$ grid that can be stacked:

\begin{lemma}\label{lem:ineq}
    If $[a]\times[b]\times[c]$ admits a box-free $2$-coloring, then 
    \[
    c \leq \frac{2\binom{a}{2}\binom{b}{2}}{\rho(a,b)}
    \]
\end{lemma}

\begin{proof}
    Suppose $[a]\times[b]\times[c]$ is 2-colored without monochromatic boxes.  Then each monochromatic red rectangle occupies a unique rectangle slot, of which there are $\binom{a}{2}\binom{b}{2}$, and the same holds for blue rectangles.  So there are $2\binom{a}{2}\binom{b}{2}$ colored rectangle slots, and $c$ layers into which $c\cdot\rho(a,b)$ such rectangles must fit, hence
    \[
    c\cdot\rho(a,b) \leq 2\binom{a}{2}\binom{b}{2}.
    \]
\end{proof}

Define \[
L(a,b) = \frac{2\binom{a}{2}\binom{b}{2}}{\rho(a,b)}.
\]

Then immediately we have $N(a,b) \leq L(a,b)$, and the colorings in Section \ref{sec:constructions} show that equality can hold.

Using the values of $\rho$ from Table \ref{tab:rho} combined with the inequality from Lemma \ref{lem:ineq}, we get some upper bounds on the largest $c$ such that $[a]\times [b]\times [c]$ is $2$-colorable. 

The advantage of this approach is that it reduces a three-dimensional extremal problem to a two-dimensional one. Rather than searching directly for box-free colorings of three-dimensional grids, we determine the minimum number $\rho(a,b)$ of monochromatic rectangles in a two-dimensional layer. A  counting argument then converts this information into an upper bound on the number of layers that can be stacked without creating a monochromatic box. In particular, this approach improves several of the upper bounds obtained by Cooper, Fenner, and Purewal, as shown in Table \ref{tab:UB}.

\begin{table}[hbt]
    \centering
    \begin{tabular}{c|c|c|c|c|c|c}
        & 3 & 4 & 5 & 6 & 7 & 8 \\
        \hline 
       3  &  &  &  &  & 126* & 84 \\ 
        \hline 
       4  &  &  &  &  & 126* & 84 \\ 
        \hline 
       5  &  &  & 100* & 75* & 52 & 46 \\ 
        \hline 
       6  &  &  &  & 75* & 52 & 46 \\ 
        \hline 
       7  &  &  &  &  & 42 (52) & 39 (45) \\ 
        \hline 
       8  &  &  &  &  &  &  39 (44) \\ 
    \end{tabular}
    \caption{Upper bounds on $N(a,b)$, the largest $c$ for which $[a]\times[b]\times[c]$  admits a box-free $2$-coloring. Parenthesized entries are the earlier bounds of Cooper, Fenner, and Purewal \cite{Fenner} that have been improved. An asterisk indicates that the displayed bound is known to be attainable by an explicit box-free coloring.}
    \label{tab:UB}
\end{table}

\section{Constructions Attaining the Layer Bound}\label{sec:constructions}

For several choices of the pair $a,b$, we have constructed colorings of $[a]\times[b]\times[c]$ such that $c$ equals the layer bound $L(a,b) = \frac{2\binom{a}{2}\binom{b}{2}}{\rho(a,b)}$. These maximal colored grids have the following dimensions:

\begin{itemize}
    \item $3\times 7 \times 126$
        \item $4\times 7 \times 126$
    \item $5\times 5 \times 100$
    \item $5\times 6 \times 75$
    \item $6\times 6 \times 75$
\end{itemize}

We will prove the  $5\times 5$ case since the proof is succinct. (Since $\rho(3,7)=1$, the $3\times 7$ case is too easy to be illustrative!) The others were handled computationally and certificates for the lower bound that have been independently checked are available at \href{https://github.com/algorithmachine/box-free-colorings}{https://github.com/algorithmachine/box-free-colorings}.

\begin{theorem}\label{initial} Any $2$-coloring of a $[5] \times [5]$ grid will contain at least two rectangles, either both of the same color or one of each.
\end{theorem}

\begin{proof}
Suppose a $5\times5$ coloring contains at most one monochromatic rectangle. As before, let $k_{c,i,j}$ count the rows in which columns $i,j$ both have color $c$, and put
\[
S=\sum_{c,i<j}k_{c,i,j}.
\]

Each row contributes $\binom r2+\binom{5-r}2$, which is an even integer at least four. Hence $S\ge20$ and $S$ is even.

Since $\binom k2\ge k-1$, the total number $R$ of monochromatic rectangles satisfies
\[
R\ge S-20.
\]

Our assumption $R\le1$ therefore forces $S=20$. Moreover, among the 20 values $k_{c,i,j}$, all but at most two must equal one.

Thus at least eight of the ten column pairs satisfy $k_{R,i,j}=k_{B,i,j}=1$. We can choose three columns whose three pairs have this property. If $r_i$ denotes the number of red entries in column $i$, then
\[
r_1+r_2=r_1+r_3=r_2+r_3=5.
\]

Summing gives $2(r_1+r_2+r_3)=15$, a contradiction.

Therefore $\rho(5,5)\ge2$.
\end{proof}


   





With this information we can now state the following theorem.

\begin{theorem}\label{thm: 5x5x100}
A $[5] \times [5]  \times [100]$ grid is $2$-colorable.
\end{theorem}

\begin{proof}
Consider the configuration below:

$$\begin{array}{|c|c|c|c|c|}
\hline
\textcolor{blue}{B} & \textcolor{blue}{B} & \textcolor{blue}{B} & \textcolor{red}{R} & \textcolor{red}{R} \\
\hline

\textcolor{blue}{B} & \textcolor{blue}{B} & \textcolor{red}{R} & \textcolor{red}{R} & \textcolor{blue}{B} \\
\hline

\textcolor{blue}{B} & \textcolor{red}{R} & \textcolor{blue}{B} & \textcolor{blue}{B} & \textcolor{red}{R} \\
\hline

\textcolor{red}{R} & \textcolor{red}{R} & \textcolor{blue}{B} & \textcolor{red}{R} & \textcolor{blue}{B} \\
\hline

\textcolor{red}{R} & \textcolor{blue}{B} & \textcolor{red}{R} & \textcolor{blue}{B} & \textcolor{red}{R} \\
\hline
\end{array}$$

It contains exactly two monochromatic rectangles, both blue.

Identify rows and columns with $\mathbb Z_5$. The two rectangles have row-pair and column-pair differences of types
\[
(1,1)\quad\text{and}\quad(2,2),
\]

where differences are considered up to sign.

Now apply all 25 independent cyclic translations of rows and columns.

Each rectangle moves through all 25 slots of its type. Consequently, these 25 layers occupy 50 distinct blue rectangle slots.

Next, multiply the column coordinates by 2 in $\mathbb Z_5$, and repeat the 25 translations.

This produces rectangles having row-pair and column-pair differences of types
\[
(1,2)\quad\text{and}\quad(2,1).
\]

These occupy the remaining 50 blue slots.

Finally, interchange the two colors in all 50 layers. This occupies all 100 red slots.

Thus we have 100 layers occupying all 200 possible colored rectangle slots exactly once. No monochromatic box occurs.
\end{proof}

\section{5 by 7 by $c$}

The constructions in the preceding section attain the layer bound exactly, suggesting the question of whether this bound is always sharp. The pair $(5,7)$ is a natural next case to investigate. Since $\rho(5,7)=8$, the layer bound gives

\[
N(5,7)\le \left\lfloor\frac{2\binom52\binom72}{8}\right\rfloor=52.
\]

We investigated this case computationally, seeking a family of 52 colorings of the $5\times7$ grid whose monochromatic rectangles occupy distinct colored rectangle slots. Our search produced a stack of 50 layers, each containing exactly eight monochromatic rectangles, with only two monochromatic boxes in the resulting three-dimensional coloring. We subsequently repaired this coloring using an exhaustive backtracking search over replacements for six of its layers, leaving the other 44 unchanged. The resulting coloring is box-free, with each layer containing exactly eight monochromatic rectangles. A certificate is available at \href{https://github.com/algorithmachine/box-free-colorings}{https://github.com/algorithmachine/box-free-colorings}.  

There are $2\binom52\binom72=420$ possible colored rectangle slots. Our construction uses exactly 400 of them, leaving just 20 unoccupied. Thus it achieves over 95\% of the theoretical packing capacity.

Combining this construction with the layer bound, we obtain

\[
50\le N(5,7)\le52.
\]

We have not found a box-free coloring with 51 or 52 layers, leading us to conjecture that $N(5,7)=50$. If true, this would demonstrate that the layer bound need not be sharp, even when every layer can attain the minimum number of monochromatic rectangles.

As $a$ and $b$ increase, the number of possible $[a]\times[b]$ layers in a box-free coloring decreases, and $[a]\times[b]\times[c]$ becomes more cube-like.  This leads us to ask what the dimensions are of the largest cube admitting a box-free coloring. 





 
\section{$n\times n\times n$ cubes}

\subsection{An upper bound}  
The inequality from Lemma \ref{lem:ineq} with $a=b=c=n$, 

 \[
    n \leq \frac{2\binom{n}{2}^2}{\rho(n,n)}
    \]

    can be rewritten 

     \[
    \rho(n,n) \leq \frac{2\binom{n}{2}^2}{n}
    \]

    to give an upper bound on $\rho(n,n)$.  If we can prove that this inequality is violated for some $n$, we immediately have that the $n$-cube is not colorable, i.e., it is $2$-guaranteed. 

    \begin{lemma}\label{lem:prob}
        We have that $\rho(a,b) \leq \frac{1}{8}\binom{a}{2}\binom{b}{2}$.
    \end{lemma}

    \begin{proof} 
We give a probabilistic argument.  (See, for example, \cite{Spencer}.) 

Color each vertex independently red or blue, each with
probability $1/2$. A given rectangle is monochromatic with
probability $2(1/2)^4=1/8$. By linearity of expectation, the
expected number of monochromatic rectangles is
\[
\frac18\binom{a}{2}\binom{b}{2}.
\]
Hence some coloring has at most this many monochromatic
rectangles, establishing the bound.
\end{proof} 

   We can use this estimate to show that our proposed method  -- that the inequality above  is violated for some $n$ --  will not work for $n<17$.   For suppose there is an $n$ such that

    \[
    \rho(n,n) > \frac{2\binom{n}{2}^2}{n}
    \]

    Then we have 

     \[
    \frac{1}{8}\binom{n}{2}^2\geq \rho(n,n) > \frac{2\binom{n}{2}^2}{n},
    \]

    and we cancel to get $n>16$. Hence we turn our attention to $n\geq 17$.

    Using a min-conflicts algorithm again, we get some upper estimates on $\rho(n,n)$. See Table \ref{tab:cube-rho}. 

    \begin{table}[hbt]
        \centering
        \begin{tabular}{c|c|c}
            $n$ &  $2\binom{n}{2}^2/n$ & empirical min-conflict upper bound on $\rho(n,n)$ \\
            \hline
             17 & 2176 & 1711 \\
             18 & 2601 & 2221 \\
             19 & 3078 & 2820\\
             20 & 3610 & 3544 \\
             21 & 4200 & 4370 \\
        \end{tabular}
        \caption{Empirical upper bounds on $\rho$}
        \label{tab:cube-rho}
    \end{table}

    $21$ is the smallest value of $n$ for which the empirical upper bound does not rule out the possibility of having 

\[
    \rho(n,n) > \frac{2\binom{n}{2}^2}{n},
    \]

    so we turn our attention to $n=21$. We will simply count rectangles, similarly to \cite{Alm}.  For fixed columns $i<j$, let $k_{c,i,j}$ denote the number of rows in which both column $i$ and column $j$ both have a dot of color $c$. 

    \begin{lemma}\label{lem:4200}
        $\sum_{c, i<j} k_{c,i,j} \geq 2100$, and the mean of the $k_{c,i,j}$s is at least $5$.
    \end{lemma}

    \begin{proof}
        Fix a row. If this row has $r$ red dots and $21-r$ blue dots, then the contribution to the sum  $\sum_{c, i<j} k_{c,i,j}$ by those dots is

        \[
        \binom{r}{2} + \binom{21-r}{2}.
        \]

        Now, 

        \[
        \binom{r}{2} + \binom{21-r}{2} = 100 + (r-10)(r-11) \geq 100. 
        \]

        Since there are $21$ rows, we get $\sum_{c, i<j} k_{c,i,j}\geq 2100$.  

        Since there are $2\binom{21}{2}=420$ terms in the sum, the average of each term is at least $2100/420 = 5$. 
    \end{proof}

    \begin{lemma}
        $\rho(21,21) \geq 4200$. 
    \end{lemma}

    \begin{proof}
        The total number of rectangles in the $21\times 21$ grid is

        \[
        \sum_{c, i<j}\binom{k_{c,i,j}}{2}.
        \]

Now, expanding the binomial coefficient around $5$ gives 
\[
\binom{k_{c,i,j}}{2} = 10 + \frac{9}{2 }(k_{c,i,j} - 5) + \frac{1}{2}(k_{c,i,j}-5)^2.
\]

Hence 

\begin{equation}\label{eq:sums}
    \sum \binom{k_{c,i,j}}{2} = \sum 10 + \frac{9}{2 }\sum (k_{c,i,j} - 5) + \frac{1}{2}\sum (k_{c,i,j}-5)^2.
\end{equation}

We note two important things: first, because the mean of the $k_{c,i,j}$s is at least $5$, $\sum (k_{c,i,j} - 5)$ is non-negative.  Second, $\sum (k_{c,i,j}-5)^2$ is non-negative because each term is.  Hence 

\begin{equation}\label{eq:sums2}
    \sum \binom{k_{c,i,j}}{2} = \sum 10 + \frac{9}{2 }\sum (k_{c,i,j} - 5) + \frac{1}{2}\sum (k_{c,i,j}-5)^2 \geq \sum 10 = 4200.
\end{equation}

    \end{proof}
    
    \begin{theorem}\label{thm:4201}
        $\rho(21,21) \geq 4201$. 
    \end{theorem}

\begin{proof}
Lemma \ref{lem:4200} established that $\rho(21,21) \geq 4200$. Suppose for contradiction that equality holds. Then there exists a coloring with exactly $4200$ rectangles. 

First, we must note that the inequality in Equation \ref{eq:sums2} holds only when every $k_{c,i,j}$ is equal to $5$.  To see this, note that because both $\sum (k_{c,i,j} - 5)$  and $\sum (k_{c,i,j}-5)^2$ are non-negative, they must be zero for equality to hold, and $\sum (k_{c,i,j}-5)^2=0$ forces all the $k_{c,i,j}$s to be $5$. 

The fact that $k_{c,i,j}=5$ for all $c$, $i$, and $j$ means that any pair of colums contains five red-red rows, five blue-blue rows, and 11 ``mixed'' rows, for a total of $21$ red dots and $21$ blue dots.  

So consider three distinct columns.  Let $r_1$, $r_2$, and $r_3$ be the number of red dots in the three colums, respectively.  By the previous paragraph, we have 

\begin{align*}
    r_1 + r_2 &= 21\\
    r_1 + r_3 &= 21\\
    r_2 + r_3 &= 21
\end{align*}

But summing gives $2(r_1+r_2+r_3) = 63$, a contradiction. 
\end{proof}

\begin{corollary}
Every $2$-coloring of the $[21]\times [21] \times [21]$ cube contains a monochromatic box.
\end{corollary}

\begin{proof}
In a box-free coloring, Lemma~4.1 gives
\[
    \rho(21,21)\le \frac{2\binom{21}{2}^{2}}{21}=4200,
\]
contradicting Theorem \ref{thm:4201}.
\end{proof}

\subsection{A lower bound on the number of monochromatic boxes}
Define a function $m:(\mathbb{Z}^+)^3 \to \mathbb{Z}^+$, where $m(a,b,c)$ is the minimum possible number of monochromatic boxes in an $a\times b \times c$ grid.  If all three dimensions are equal, we'll write $m(n)$ for $m(n,n,n)$.  

\begin{lemma}\label{lem:m}
    \[
m(a,b,c)\ge c\rho(a,b)-2\binom a2\binom b2.
\]

In particular, for cubes we have 

  \[
m(n)\geq n\rho(n,n)-2\binom n2^2.
\]
\end{lemma}

\begin{proof}
    Consider the $c$ layers of $[a]\times [b]$ grids.

Each layer contains at least $\rho(a,b)$ monochromatic rectangles, giving at least $c\rho(a,b)$ colored rectangles altogether, counted with multiplicity.

But there are only $2\binom a2\binom b2$ possible colored rectangles: two choices of color and $\binom a2\binom b2$ choices of slot.

Each occurrence beyond the first of a particular colored rectangle creates at least one monochromatic box with an earlier occurrence. Consequently,
\[
m(a,b,c)\geq c\rho(a,b)-2\binom a2\binom b2.
\]
\end{proof}

\begin{corollary}
    \[
m(21)\geq 21.
\]
\end{corollary}

\begin{proof}
    We have established that $\rho(21,21) \geq 4201.$ Hence by Lemma \ref{lem:m} we have 

      \[
m(21)\geq 21\rho(21,21)-2\binom {21}{2}^2 \geq 21(4201) - 88200 = 88221-88200 = 21.
\]
\end{proof}

\subsection{$11\times 11\times 11$}
We spent considerable effort trying to determine whether the $11\times 11\times 11$ cube was colorable, and in this section we detail our efforts.

We investigated the $11\times 11\times 11$ cube by computational local
search, with the goal of minimizing the number of monochromatic boxes.
A coloring has $11^3=1331$ vertices, and there are
$\binom{11}{2}^3=166375$ possible axis-parallel boxes.  Since each vertex
belongs to $10^3=1000$ boxes, the effect of a single vertex flip can be
computed by inspecting only those $1000$ boxes.

Our  approach was \emph{min-conflicts}: we tried to minimize the ``conflicts'', i.e., bad boxes.
Beginning with a random two-coloring, we chose a monochromatic box uniformly
at random and examined the eight color changes obtained by flipping one of
its corners.  We performed a flip that produced the fewest total
number of monochromatic boxes, breaking ties randomly.  We also incorporated
perturbations to escape local minima: 
a run was perturbed by $50$ randomly chosen vertex flips after $2000$
steps if no improvement was seen.  

\begin{figure}
    \centering
    \includegraphics[width=1.1\linewidth]{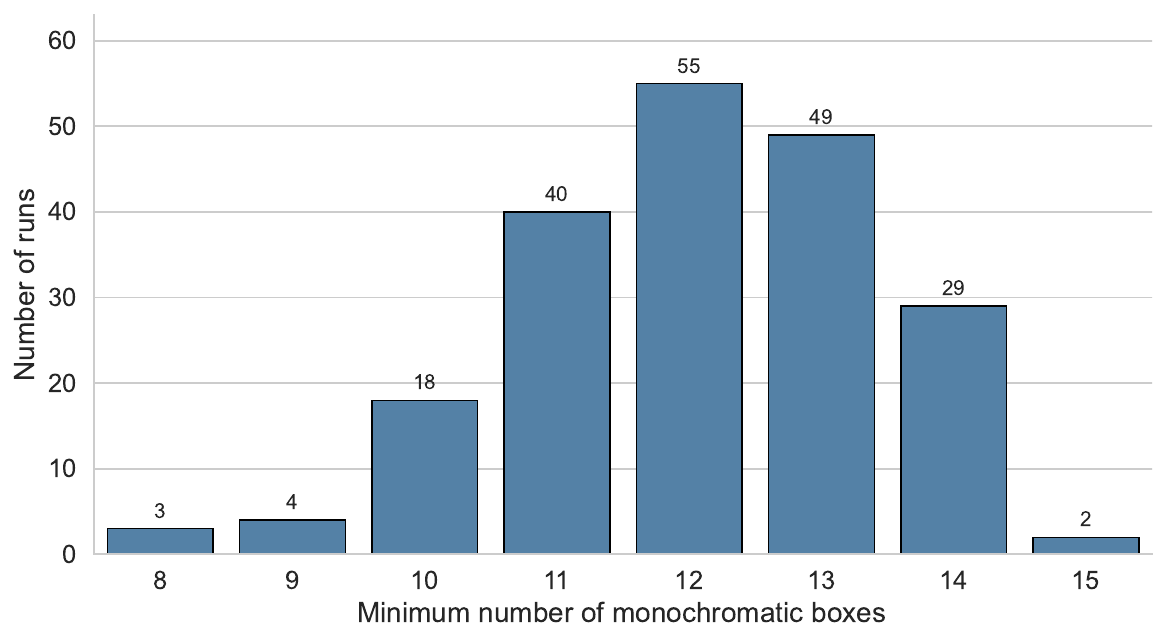}
    \caption{Smallest number of monochromatic boxes in $200$ independent $500,000$-step searches
for two-colorings of the $11\times11\times11$ cube.}
    \label{fig:11cube-search}
\end{figure}

In one experiment, we performed $200$ independent runs, allowing $500,000$
search steps per run.  Figure~\ref{fig:11cube-search} gives the smallest
number of monochromatic boxes attained in each run.

The longer searches were significant: the three eight-box colorings in this
batch first appeared at steps $169319$, $193279$, and $432911$, and all four
nine-box colorings appeared after step $300000$.  Across our experiments, we
found four independently obtained colorings with precisely eight
monochromatic boxes.  These colorings were checked to be inequivalent under
permutations of the coordinate axes, independent permutations of the
coordinate values, and interchange of the two colors.  They had $669$, $665$,
$661$, and $667$ vertices of color~1, respectively.  In two of these
colorings, the eight monochromatic boxes were pairwise vertex-disjoint.

Consequently,
\[
    m(11)\leq 8.
\]
  We
therefore conjecture that $m(11)=8$, which would imply that
$10\times10\times10$ is the largest box-free cube.

\section{Open problems}

The constructions and computations in this paper suggest several directions for further work. We record six problems that seem particularly likely to lead either to sharper theoretical results or to tractable computational projects.

\begin{problem}
Determine $m(11)$.
\end{problem}

Recall that $m(n)$ denotes the minimum number of monochromatic boxes in a $2$-coloring of $[n]^3$. We have exhibited several colorings of $[11]^3$ containing exactly eight monochromatic boxes, and hence $m(11)\leq 8$.

We found no coloring with fewer than eight, leading us to conjecture that  $m(11) = 8$.

A proof would, in particular, show that $[10]^3$ is the largest cube admitting a box-free $2$-coloring.  

\begin{problem}
Determine the behavior of $m(n)$ asymptotically as $n\to\infty$.
\end{problem}

Even after $m(11)$, the values $m(12),m(13),\ldots$ are unknown. More broadly, essentially nothing is known about the order of growth of $m(n)$. Since $[n]^3$ contains $\binom{n}{2}^3$ boxes, one may ask whether

$$
\frac{m(n)}{\binom {n}{2}^3}
$$

has a limit, or even whether $m(n)$ has polynomial order $n^\alpha$ for some $\alpha$. Lemma \ref{lem:prob} gives an upper bound on $\rho$, but it is based on random colorings, so is surely too large. 

\begin{problem}
For fixed $a$ and $b$, determine the largest $c$ for which $[a]\times[b]\times[c]$ admits a box-free $2$-coloring.
\end{problem}

Define

$$
N(a,b)=\max\{c:[a]\times[b]\times[c]\text{ admits a box-free $2$-coloring}\}.
$$

The layer-counting argument gives

$$
N(a,b)\leq\frac{2\binom a2\binom b2}{\rho(a,b)},
$$

where $\rho(a,b)$ is the minimum number of monochromatic rectangles in a $2$-coloring of an $a\times b$ grid. Our constructions attain this bound for several pairs, including $(3,7)$, $(4,7)$, $(5,5)$, $(5,6)$, and $(6,6)$. It would be interesting to determine $N(a,b)$ systematically and, in particular, to find the smallest pair $(a,b)$ for which this upper bound -- the layer bound -- is not attainable. The pair $(5,7)$ is a natural first candidate. We know $50 \leq N(5,7) \leq 52$.  

\begin{problem}
Determine $\rho(a,b)$, study its asymptotic behavior, and characterize the pairs for which its convexity lower bound is sharp.
\end{problem}

The function $\rho(a,b)$ is an extremal problem of independent interest. Even the diagonal sequence

$$
\rho(n,n)=0,0,0,0,2,6,21,40,\ldots
$$

quickly becomes nontrivial. The general convexity argument gives an explicit lower bound for $\rho(a,b)$, and this bound is attained by every nonzero value in our Table \ref{tab:rho} except $\rho(5,5)=2$. It is therefore natural to ask when the bound is sharp and to characterize its equality cases.

\begin{problem}
Strengthen the layer-counting method by incorporating the structure of monochromatic rectangles within and between layers.
\end{problem}

Every layer of a coloring of $[a]\times[b]\times[c]$ contains at least $\rho(a,b)$ monochromatic rectangles, while any fixed rectangle-position and color can occur in at most one layer of a box-free coloring. The resulting upper bound on $N(a,b)$ uses only the total number of rectangles in each layer and discards most of their structure.

A sharper argument might account for the distribution of rectangles among row-pairs, column-pairs, and colors, or for restrictions on which minimal layers can coexist.  One particularly attractive goal is  to prove that $[11]^3$ is not box-free, or even to obtain a meaningful lower bound on $m(11)$.

\begin{problem}
Table \ref{tab:fenner-upper-bounds} from \cite{Fenner} has more entries than our Table \ref{tab:UB}, although some entries in the latter improve upon those in the former via use of the layer bound.  Computing for values of $\rho$ would likely lead to more such improvements.    
\end{problem}

\section*{AI Acknowledgement}

ChatGPT (OpenAI, GPT-6) was used by the  authors for coding assistance, idea generation, computational exploration, figure creation, and copyediting. 

The authors have checked terms of use for ChatGPT, and therefore confirm suitability for publication.  Originality and accuracy of content has been confirmed by the authors. The authors take full responsibility for the integrity of the whole content, including accuracy of references.  All claims not proven in the document have available certificates that can be verified independently. 

\section*{Data Availability}

All data, in the form of certificates, have been made available at 
\begin{center}
    \href{https://github.com/algorithmachine/box-free-colorings}{https://github.com/algorithmachine/box-free-colorings}.  
\end{center}


\appendix

\section{$\rho$ certificates}\label{appendix:rho-certs}

\label{app:rho-certificates}
For each nonzero entry in Table~2, the following binary matrix gives a
coloring attaining the stated value of $\rho(a,b)$, with black representing
color~$1$ and white representing color~$0$. Rows are displayed from top to
bottom and columns from left to right. The adjacent binary strings provide
an exact, machine-readable transcription. Each claimed rectangle count was
checked by enumerating all $\binom{a}{2}\binom{b}{2}$ possible rectangles. These images were generated by ChatGPT. 

\medskip

\begin{minipage}[t]{0.30\linewidth}
\centering
\vbox{\hbox to \linewidth{\hfil $\boldsymbol{\rho(3,7)=1}$ \hfil}}
\vspace{0.45em}
\begin{tabular}{@{}c@{\hspace{0.5em}}l@{}}
\begin{tikzpicture}[x=0.275cm,y=0.275cm,baseline=(current bounding box.center)]
\fill[black] (1,0) rectangle ++(1,-1);
\fill[black] (3,0) rectangle ++(1,-1);
\fill[black] (5,0) rectangle ++(1,-1);
\fill[black] (0,-1) rectangle ++(1,-1);
\fill[black] (3,-1) rectangle ++(1,-1);
\fill[black] (4,-1) rectangle ++(1,-1);
\fill[black] (6,-1) rectangle ++(1,-1);
\fill[black] (2,-2) rectangle ++(1,-1);
\fill[black] (4,-2) rectangle ++(1,-1);
\fill[black] (5,-2) rectangle ++(1,-1);
\fill[black] (6,-2) rectangle ++(1,-1);
\draw[gray!65,line width=0.24pt,step=1] (0,-3) grid (7,0);
\end{tikzpicture}
 & 
{\ttfamily\scriptsize\begin{tabular}{@{}l@{}}
0101010 \\
1001101 \\
0010111
\end{tabular}}
\end{tabular}
\end{minipage}
\hfill
\begin{minipage}[t]{0.30\linewidth}
\centering
\vbox{\hbox to \linewidth{\hfil $\boldsymbol{\rho(3,8)=2}$ \hfil}}
\vspace{0.45em}
\begin{tabular}{@{}c@{\hspace{0.5em}}l@{}}
\begin{tikzpicture}[x=0.275cm,y=0.275cm,baseline=(current bounding box.center)]
\fill[black] (0,0) rectangle ++(1,-1);
\fill[black] (1,0) rectangle ++(1,-1);
\fill[black] (6,0) rectangle ++(1,-1);
\fill[black] (7,0) rectangle ++(1,-1);
\fill[black] (1,-1) rectangle ++(1,-1);
\fill[black] (3,-1) rectangle ++(1,-1);
\fill[black] (5,-1) rectangle ++(1,-1);
\fill[black] (2,-2) rectangle ++(1,-1);
\fill[black] (4,-2) rectangle ++(1,-1);
\fill[black] (5,-2) rectangle ++(1,-1);
\fill[black] (7,-2) rectangle ++(1,-1);
\draw[gray!65,line width=0.24pt,step=1] (0,-3) grid (8,0);
\end{tikzpicture}
 & 
{\ttfamily\scriptsize\begin{tabular}{@{}l@{}}
11000011 \\
01010100 \\
00101101
\end{tabular}}
\end{tabular}
\end{minipage}
\hfill
\begin{minipage}[t]{0.30\linewidth}
\centering
\vbox{\hbox to \linewidth{\hfil $\boldsymbol{\rho(4,7)=2}$ \hfil}}
\vspace{0.45em}
\begin{tabular}{@{}c@{\hspace{0.5em}}l@{}}
\begin{tikzpicture}[x=0.275cm,y=0.275cm,baseline=(current bounding box.center)]
\fill[black] (2,0) rectangle ++(1,-1);
\fill[black] (3,0) rectangle ++(1,-1);
\fill[black] (4,0) rectangle ++(1,-1);
\fill[black] (2,-1) rectangle ++(1,-1);
\fill[black] (5,-1) rectangle ++(1,-1);
\fill[black] (6,-1) rectangle ++(1,-1);
\fill[black] (0,-2) rectangle ++(1,-1);
\fill[black] (1,-2) rectangle ++(1,-1);
\fill[black] (3,-2) rectangle ++(1,-1);
\fill[black] (6,-2) rectangle ++(1,-1);
\fill[black] (0,-3) rectangle ++(1,-1);
\fill[black] (1,-3) rectangle ++(1,-1);
\fill[black] (4,-3) rectangle ++(1,-1);
\fill[black] (5,-3) rectangle ++(1,-1);
\draw[gray!65,line width=0.24pt,step=1] (0,-4) grid (7,0);
\end{tikzpicture}
 & 
{\ttfamily\scriptsize\begin{tabular}{@{}l@{}}
0011100 \\
0010011 \\
1101001 \\
1100110
\end{tabular}}
\end{tabular}
\end{minipage}
\par\vspace{1.15em}
\begin{minipage}[t]{0.30\linewidth}
\centering
\vbox{\hbox to \linewidth{\hfil $\boldsymbol{\rho(4,8)=4}$ \hfil}}
\vspace{0.45em}
\begin{tabular}{@{}c@{\hspace{0.5em}}l@{}}
\begin{tikzpicture}[x=0.275cm,y=0.275cm,baseline=(current bounding box.center)]
\fill[black] (0,0) rectangle ++(1,-1);
\fill[black] (5,0) rectangle ++(1,-1);
\fill[black] (6,0) rectangle ++(1,-1);
\fill[black] (7,0) rectangle ++(1,-1);
\fill[black] (0,-1) rectangle ++(1,-1);
\fill[black] (1,-1) rectangle ++(1,-1);
\fill[black] (2,-1) rectangle ++(1,-1);
\fill[black] (3,-1) rectangle ++(1,-1);
\fill[black] (7,-1) rectangle ++(1,-1);
\fill[black] (2,-2) rectangle ++(1,-1);
\fill[black] (4,-2) rectangle ++(1,-1);
\fill[black] (5,-2) rectangle ++(1,-1);
\fill[black] (1,-3) rectangle ++(1,-1);
\fill[black] (3,-3) rectangle ++(1,-1);
\fill[black] (4,-3) rectangle ++(1,-1);
\fill[black] (6,-3) rectangle ++(1,-1);
\draw[gray!65,line width=0.24pt,step=1] (0,-4) grid (8,0);
\end{tikzpicture}
 & 
{\ttfamily\scriptsize\begin{tabular}{@{}l@{}}
10000111 \\
11110001 \\
00101100 \\
01011010
\end{tabular}}
\end{tabular}
\end{minipage}
\hfill
\begin{minipage}[t]{0.30\linewidth}
\centering
\vbox{\hbox to \linewidth{\hfil $\boldsymbol{\rho(5,5)=2}$ \hfil}}
\vspace{0.45em}
\begin{tabular}{@{}c@{\hspace{0.5em}}l@{}}
\begin{tikzpicture}[x=0.275cm,y=0.275cm,baseline=(current bounding box.center)]
\fill[black] (0,0) rectangle ++(1,-1);
\fill[black] (4,0) rectangle ++(1,-1);
\fill[black] (0,-1) rectangle ++(1,-1);
\fill[black] (1,-1) rectangle ++(1,-1);
\fill[black] (1,-2) rectangle ++(1,-1);
\fill[black] (3,-2) rectangle ++(1,-1);
\fill[black] (4,-2) rectangle ++(1,-1);
\fill[black] (1,-3) rectangle ++(1,-1);
\fill[black] (2,-3) rectangle ++(1,-1);
\fill[black] (0,-4) rectangle ++(1,-1);
\fill[black] (2,-4) rectangle ++(1,-1);
\fill[black] (3,-4) rectangle ++(1,-1);
\draw[gray!65,line width=0.24pt,step=1] (0,-5) grid (5,0);
\end{tikzpicture}
 & 
{\ttfamily\scriptsize\begin{tabular}{@{}l@{}}
10001 \\
11000 \\
01011 \\
01100 \\
10110
\end{tabular}}
\end{tabular}
\end{minipage}
\hfill
\begin{minipage}[t]{0.30\linewidth}
\centering
\vbox{\hbox to \linewidth{\hfil $\boldsymbol{\rho(5,6)=4}$ \hfil}}
\vspace{0.45em}
\begin{tabular}{@{}c@{\hspace{0.5em}}l@{}}
\begin{tikzpicture}[x=0.275cm,y=0.275cm,baseline=(current bounding box.center)]
\fill[black] (0,0) rectangle ++(1,-1);
\fill[black] (3,0) rectangle ++(1,-1);
\fill[black] (5,0) rectangle ++(1,-1);
\fill[black] (1,-1) rectangle ++(1,-1);
\fill[black] (4,-1) rectangle ++(1,-1);
\fill[black] (5,-1) rectangle ++(1,-1);
\fill[black] (0,-2) rectangle ++(1,-1);
\fill[black] (1,-2) rectangle ++(1,-1);
\fill[black] (3,-2) rectangle ++(1,-1);
\fill[black] (4,-2) rectangle ++(1,-1);
\fill[black] (0,-3) rectangle ++(1,-1);
\fill[black] (2,-3) rectangle ++(1,-1);
\fill[black] (4,-3) rectangle ++(1,-1);
\fill[black] (1,-4) rectangle ++(1,-1);
\fill[black] (2,-4) rectangle ++(1,-1);
\fill[black] (3,-4) rectangle ++(1,-1);
\draw[gray!65,line width=0.24pt,step=1] (0,-5) grid (6,0);
\end{tikzpicture}
 & 
{\ttfamily\scriptsize\begin{tabular}{@{}l@{}}
100101 \\
010011 \\
110110 \\
101010 \\
011100
\end{tabular}}
\end{tabular}
\end{minipage}
\par\vspace{1.15em}
\begin{minipage}[t]{0.30\linewidth}
\centering
\vbox{\hbox to \linewidth{\hfil $\boldsymbol{\rho(5,7)=8}$ \hfil}}
\vspace{0.45em}
\begin{tabular}{@{}c@{\hspace{0.5em}}l@{}}
\begin{tikzpicture}[x=0.275cm,y=0.275cm,baseline=(current bounding box.center)]
\fill[black] (1,0) rectangle ++(1,-1);
\fill[black] (2,0) rectangle ++(1,-1);
\fill[black] (5,0) rectangle ++(1,-1);
\fill[black] (6,0) rectangle ++(1,-1);
\fill[black] (0,-1) rectangle ++(1,-1);
\fill[black] (5,-1) rectangle ++(1,-1);
\fill[black] (6,-1) rectangle ++(1,-1);
\fill[black] (2,-2) rectangle ++(1,-1);
\fill[black] (3,-2) rectangle ++(1,-1);
\fill[black] (4,-2) rectangle ++(1,-1);
\fill[black] (6,-2) rectangle ++(1,-1);
\fill[black] (0,-3) rectangle ++(1,-1);
\fill[black] (1,-3) rectangle ++(1,-1);
\fill[black] (4,-3) rectangle ++(1,-1);
\fill[black] (0,-4) rectangle ++(1,-1);
\fill[black] (2,-4) rectangle ++(1,-1);
\fill[black] (3,-4) rectangle ++(1,-1);
\draw[gray!65,line width=0.24pt,step=1] (0,-5) grid (7,0);
\end{tikzpicture}
 & 
{\ttfamily\scriptsize\begin{tabular}{@{}l@{}}
0110011 \\
1000011 \\
0011101 \\
1100100 \\
1011000
\end{tabular}}
\end{tabular}
\end{minipage}
\hfill
\begin{minipage}[t]{0.30\linewidth}
\centering
\vbox{\hbox to \linewidth{\hfil $\boldsymbol{\rho(5,8)=12}$ \hfil}}
\vspace{0.45em}
\begin{tabular}{@{}c@{\hspace{0.5em}}l@{}}
\begin{tikzpicture}[x=0.275cm,y=0.275cm,baseline=(current bounding box.center)]
\fill[black] (0,0) rectangle ++(1,-1);
\fill[black] (4,0) rectangle ++(1,-1);
\fill[black] (5,0) rectangle ++(1,-1);
\fill[black] (7,0) rectangle ++(1,-1);
\fill[black] (0,-1) rectangle ++(1,-1);
\fill[black] (2,-1) rectangle ++(1,-1);
\fill[black] (3,-1) rectangle ++(1,-1);
\fill[black] (6,-1) rectangle ++(1,-1);
\fill[black] (7,-1) rectangle ++(1,-1);
\fill[black] (0,-2) rectangle ++(1,-1);
\fill[black] (1,-2) rectangle ++(1,-1);
\fill[black] (2,-2) rectangle ++(1,-1);
\fill[black] (1,-3) rectangle ++(1,-1);
\fill[black] (3,-3) rectangle ++(1,-1);
\fill[black] (5,-3) rectangle ++(1,-1);
\fill[black] (7,-3) rectangle ++(1,-1);
\fill[black] (1,-4) rectangle ++(1,-1);
\fill[black] (3,-4) rectangle ++(1,-1);
\fill[black] (4,-4) rectangle ++(1,-1);
\fill[black] (6,-4) rectangle ++(1,-1);
\draw[gray!65,line width=0.24pt,step=1] (0,-5) grid (8,0);
\end{tikzpicture}
 & 
{\ttfamily\scriptsize\begin{tabular}{@{}l@{}}
10001101 \\
10110011 \\
11100000 \\
01010101 \\
01011010
\end{tabular}}
\end{tabular}
\end{minipage}
\hfill
\begin{minipage}[t]{0.30\linewidth}
\centering
\vbox{\hbox to \linewidth{\hfil $\boldsymbol{\rho(6,6)=6}$ \hfil}}
\vspace{0.45em}
\begin{tabular}{@{}c@{\hspace{0.5em}}l@{}}
\begin{tikzpicture}[x=0.275cm,y=0.275cm,baseline=(current bounding box.center)]
\fill[black] (0,0) rectangle ++(1,-1);
\fill[black] (1,0) rectangle ++(1,-1);
\fill[black] (5,0) rectangle ++(1,-1);
\fill[black] (0,-1) rectangle ++(1,-1);
\fill[black] (2,-1) rectangle ++(1,-1);
\fill[black] (3,-1) rectangle ++(1,-1);
\fill[black] (1,-2) rectangle ++(1,-1);
\fill[black] (2,-2) rectangle ++(1,-1);
\fill[black] (4,-2) rectangle ++(1,-1);
\fill[black] (2,-3) rectangle ++(1,-1);
\fill[black] (3,-3) rectangle ++(1,-1);
\fill[black] (5,-3) rectangle ++(1,-1);
\fill[black] (0,-4) rectangle ++(1,-1);
\fill[black] (4,-4) rectangle ++(1,-1);
\fill[black] (5,-4) rectangle ++(1,-1);
\fill[black] (1,-5) rectangle ++(1,-1);
\fill[black] (3,-5) rectangle ++(1,-1);
\fill[black] (4,-5) rectangle ++(1,-1);
\draw[gray!65,line width=0.24pt,step=1] (0,-6) grid (6,0);
\end{tikzpicture}
 & 
{\ttfamily\scriptsize\begin{tabular}{@{}l@{}}
110001 \\
101100 \\
011010 \\
001101 \\
100011 \\
010110
\end{tabular}}
\end{tabular}
\end{minipage}
\par\vspace{1.15em}
\begin{minipage}[t]{0.30\linewidth}
\centering
\vbox{\hbox to \linewidth{\hfil $\boldsymbol{\rho(6,7)=12}$ \hfil}}
\vspace{0.45em}
\begin{tabular}{@{}c@{\hspace{0.5em}}l@{}}
\begin{tikzpicture}[x=0.275cm,y=0.275cm,baseline=(current bounding box.center)]
\fill[black] (2,0) rectangle ++(1,-1);
\fill[black] (3,0) rectangle ++(1,-1);
\fill[black] (4,0) rectangle ++(1,-1);
\fill[black] (5,0) rectangle ++(1,-1);
\fill[black] (0,-1) rectangle ++(1,-1);
\fill[black] (2,-1) rectangle ++(1,-1);
\fill[black] (5,-1) rectangle ++(1,-1);
\fill[black] (6,-1) rectangle ++(1,-1);
\fill[black] (0,-2) rectangle ++(1,-1);
\fill[black] (3,-2) rectangle ++(1,-1);
\fill[black] (4,-2) rectangle ++(1,-1);
\fill[black] (1,-3) rectangle ++(1,-1);
\fill[black] (3,-3) rectangle ++(1,-1);
\fill[black] (6,-3) rectangle ++(1,-1);
\fill[black] (1,-4) rectangle ++(1,-1);
\fill[black] (4,-4) rectangle ++(1,-1);
\fill[black] (5,-4) rectangle ++(1,-1);
\fill[black] (6,-4) rectangle ++(1,-1);
\fill[black] (0,-5) rectangle ++(1,-1);
\fill[black] (1,-5) rectangle ++(1,-1);
\fill[black] (2,-5) rectangle ++(1,-1);
\draw[gray!65,line width=0.24pt,step=1] (0,-6) grid (7,0);
\end{tikzpicture}
 & 
{\ttfamily\scriptsize\begin{tabular}{@{}l@{}}
0011110 \\
1010011 \\
1001100 \\
0101001 \\
0100111 \\
1110000
\end{tabular}}
\end{tabular}
\end{minipage}
\hfill
\begin{minipage}[t]{0.30\linewidth}
\centering
\vbox{\hbox to \linewidth{\hfil $\boldsymbol{\rho(6,8)=18}$ \hfil}}
\vspace{0.45em}
\begin{tabular}{@{}c@{\hspace{0.5em}}l@{}}
\begin{tikzpicture}[x=0.275cm,y=0.275cm,baseline=(current bounding box.center)]
\fill[black] (2,0) rectangle ++(1,-1);
\fill[black] (4,0) rectangle ++(1,-1);
\fill[black] (5,0) rectangle ++(1,-1);
\fill[black] (7,0) rectangle ++(1,-1);
\fill[black] (0,-1) rectangle ++(1,-1);
\fill[black] (1,-1) rectangle ++(1,-1);
\fill[black] (3,-1) rectangle ++(1,-1);
\fill[black] (5,-1) rectangle ++(1,-1);
\fill[black] (0,-2) rectangle ++(1,-1);
\fill[black] (4,-2) rectangle ++(1,-1);
\fill[black] (5,-2) rectangle ++(1,-1);
\fill[black] (6,-2) rectangle ++(1,-1);
\fill[black] (1,-3) rectangle ++(1,-1);
\fill[black] (3,-3) rectangle ++(1,-1);
\fill[black] (4,-3) rectangle ++(1,-1);
\fill[black] (6,-3) rectangle ++(1,-1);
\fill[black] (7,-3) rectangle ++(1,-1);
\fill[black] (0,-4) rectangle ++(1,-1);
\fill[black] (1,-4) rectangle ++(1,-1);
\fill[black] (2,-4) rectangle ++(1,-1);
\fill[black] (7,-4) rectangle ++(1,-1);
\fill[black] (2,-5) rectangle ++(1,-1);
\fill[black] (3,-5) rectangle ++(1,-1);
\fill[black] (6,-5) rectangle ++(1,-1);
\draw[gray!65,line width=0.24pt,step=1] (0,-6) grid (8,0);
\end{tikzpicture}
 & 
{\ttfamily\scriptsize\begin{tabular}{@{}l@{}}
00101101 \\
11010100 \\
10001110 \\
01011011 \\
11100001 \\
00110010
\end{tabular}}
\end{tabular}
\end{minipage}
\hfill
\begin{minipage}[t]{0.30\linewidth}
\centering
\vbox{\hbox to \linewidth{\hfil $\boldsymbol{\rho(7,7)=21}$ \hfil}}
\vspace{0.45em}
\begin{tabular}{@{}c@{\hspace{0.5em}}l@{}}
\begin{tikzpicture}[x=0.275cm,y=0.275cm,baseline=(current bounding box.center)]
\fill[black] (0,0) rectangle ++(1,-1);
\fill[black] (1,0) rectangle ++(1,-1);
\fill[black] (3,0) rectangle ++(1,-1);
\fill[black] (1,-1) rectangle ++(1,-1);
\fill[black] (2,-1) rectangle ++(1,-1);
\fill[black] (4,-1) rectangle ++(1,-1);
\fill[black] (2,-2) rectangle ++(1,-1);
\fill[black] (3,-2) rectangle ++(1,-1);
\fill[black] (5,-2) rectangle ++(1,-1);
\fill[black] (3,-3) rectangle ++(1,-1);
\fill[black] (4,-3) rectangle ++(1,-1);
\fill[black] (6,-3) rectangle ++(1,-1);
\fill[black] (0,-4) rectangle ++(1,-1);
\fill[black] (4,-4) rectangle ++(1,-1);
\fill[black] (5,-4) rectangle ++(1,-1);
\fill[black] (1,-5) rectangle ++(1,-1);
\fill[black] (5,-5) rectangle ++(1,-1);
\fill[black] (6,-5) rectangle ++(1,-1);
\fill[black] (0,-6) rectangle ++(1,-1);
\fill[black] (2,-6) rectangle ++(1,-1);
\fill[black] (6,-6) rectangle ++(1,-1);
\draw[gray!65,line width=0.24pt,step=1] (0,-7) grid (7,0);
\end{tikzpicture}
 & 
{\ttfamily\scriptsize\begin{tabular}{@{}l@{}}
1101000 \\
0110100 \\
0011010 \\
0001101 \\
1000110 \\
0100011 \\
1010001
\end{tabular}}
\end{tabular}
\end{minipage}
\par\vspace{1.15em}
\begin{minipage}[t]{0.30\linewidth}
\centering
\vbox{\hbox to \linewidth{\hfil $\boldsymbol{\rho(7,8)=30}$ \hfil}}
\vspace{0.45em}
\begin{tabular}{@{}c@{\hspace{0.5em}}l@{}}
\begin{tikzpicture}[x=0.275cm,y=0.275cm,baseline=(current bounding box.center)]
\fill[black] (0,0) rectangle ++(1,-1);
\fill[black] (1,0) rectangle ++(1,-1);
\fill[black] (2,0) rectangle ++(1,-1);
\fill[black] (6,0) rectangle ++(1,-1);
\fill[black] (0,-1) rectangle ++(1,-1);
\fill[black] (2,-1) rectangle ++(1,-1);
\fill[black] (5,-1) rectangle ++(1,-1);
\fill[black] (7,-1) rectangle ++(1,-1);
\fill[black] (1,-2) rectangle ++(1,-1);
\fill[black] (3,-2) rectangle ++(1,-1);
\fill[black] (5,-2) rectangle ++(1,-1);
\fill[black] (7,-2) rectangle ++(1,-1);
\fill[black] (0,-3) rectangle ++(1,-1);
\fill[black] (4,-3) rectangle ++(1,-1);
\fill[black] (5,-3) rectangle ++(1,-1);
\fill[black] (6,-3) rectangle ++(1,-1);
\fill[black] (2,-4) rectangle ++(1,-1);
\fill[black] (3,-4) rectangle ++(1,-1);
\fill[black] (6,-4) rectangle ++(1,-1);
\fill[black] (7,-4) rectangle ++(1,-1);
\fill[black] (0,-5) rectangle ++(1,-1);
\fill[black] (3,-5) rectangle ++(1,-1);
\fill[black] (4,-5) rectangle ++(1,-1);
\fill[black] (7,-5) rectangle ++(1,-1);
\fill[black] (1,-6) rectangle ++(1,-1);
\fill[black] (2,-6) rectangle ++(1,-1);
\fill[black] (3,-6) rectangle ++(1,-1);
\fill[black] (4,-6) rectangle ++(1,-1);
\draw[gray!65,line width=0.24pt,step=1] (0,-7) grid (8,0);
\end{tikzpicture}
 & 
{\ttfamily\scriptsize\begin{tabular}{@{}l@{}}
11100010 \\
10100101 \\
01010101 \\
10001110 \\
00110011 \\
10011001 \\
01111000
\end{tabular}}
\end{tabular}
\end{minipage}
\hfill
\begin{minipage}[t]{0.30\linewidth}
\centering
\vbox{\hbox to \linewidth{\hfil $\boldsymbol{\rho(8,8)=40}$ \hfil}}
\vspace{0.45em}
\begin{tabular}{@{}c@{\hspace{0.5em}}l@{}}
\begin{tikzpicture}[x=0.275cm,y=0.275cm,baseline=(current bounding box.center)]
\fill[black] (0,0) rectangle ++(1,-1);
\fill[black] (1,0) rectangle ++(1,-1);
\fill[black] (2,0) rectangle ++(1,-1);
\fill[black] (4,0) rectangle ++(1,-1);
\fill[black] (1,-1) rectangle ++(1,-1);
\fill[black] (2,-1) rectangle ++(1,-1);
\fill[black] (3,-1) rectangle ++(1,-1);
\fill[black] (5,-1) rectangle ++(1,-1);
\fill[black] (2,-2) rectangle ++(1,-1);
\fill[black] (3,-2) rectangle ++(1,-1);
\fill[black] (4,-2) rectangle ++(1,-1);
\fill[black] (6,-2) rectangle ++(1,-1);
\fill[black] (3,-3) rectangle ++(1,-1);
\fill[black] (4,-3) rectangle ++(1,-1);
\fill[black] (5,-3) rectangle ++(1,-1);
\fill[black] (7,-3) rectangle ++(1,-1);
\fill[black] (0,-4) rectangle ++(1,-1);
\fill[black] (4,-4) rectangle ++(1,-1);
\fill[black] (5,-4) rectangle ++(1,-1);
\fill[black] (6,-4) rectangle ++(1,-1);
\fill[black] (1,-5) rectangle ++(1,-1);
\fill[black] (5,-5) rectangle ++(1,-1);
\fill[black] (6,-5) rectangle ++(1,-1);
\fill[black] (7,-5) rectangle ++(1,-1);
\fill[black] (0,-6) rectangle ++(1,-1);
\fill[black] (2,-6) rectangle ++(1,-1);
\fill[black] (6,-6) rectangle ++(1,-1);
\fill[black] (7,-6) rectangle ++(1,-1);
\fill[black] (0,-7) rectangle ++(1,-1);
\fill[black] (1,-7) rectangle ++(1,-1);
\fill[black] (3,-7) rectangle ++(1,-1);
\fill[black] (7,-7) rectangle ++(1,-1);
\draw[gray!65,line width=0.24pt,step=1] (0,-8) grid (8,0);
\end{tikzpicture}
 & 
{\ttfamily\scriptsize\begin{tabular}{@{}l@{}}
11101000 \\
01110100 \\
00111010 \\
00011101 \\
10001110 \\
01000111 \\
10100011 \\
11010001
\end{tabular}}
\end{tabular}
\end{minipage}
\hfill
\par\medskip

\bibliographystyle{apalike}

\bibliography{refs}

\end{document}